# QUADRATIC DISTANCES ON PROBABILITIES: A UNIFIED FOUNDATION

By Bruce G. Lindsay,[1] Marianthi Markatou,[2] Surajit Ray, Ke Yang and Shu-Chuan Chen

*Pennsylvania State University, Columbia University, Boston University, Cytokinetics Inc. and National Cheng Kung University*

This work builds a unified framework for the study of quadratic form distance measures as they are used in assessing the goodness of fit of models. Many important procedures have this structure, but the theory for these methods is dispersed and incomplete. Central to the statistical analysis of these distances is the spectral decomposition of the kernel that generates the distance. We show how this determines the limiting distribution of natural goodness-of-fit tests. Additionally, we develop a new notion, the spectral degrees of freedom of the test, based on this decomposition. The degrees of freedom are easy to compute and estimate, and can be used as a guide in the construction of useful procedures in this class.

**1. Introduction.** Modern scientific work has presented statistics with many important challenges, but of particular importance are the challenges presented by "large magnitude," both in the dimension of data vectors and in the number of vectors [see Lindsay, Kettenring and Siegmund (2004)]. Assessment of the fit of a model in such a situation can be challenging.

Model fit assessment is usually based, either explicitly or implicitly, on measures of distance $d(F, G)$ between probability measures $F$ and $G$. Our foundation stones will be quadratic distance measures. This class is characterized by the simple quadratic form structure

$$d_K(F, G) = \iint K_G(s, t) \, d(F - G)(s) \, d(F - G)(t),$$

Received March 2006; revised April 2007.
[1]Supported by NSF Grant DMS-04-05637.
[2]Supported by NSF Grant DMS-05-04957.

that is adaptable through the choice of a nonnegative definite kernel $K_G(s,t)$. This form is asymmetric in $F$ and $G$; here $G$ will often be a distribution whose goodness of fit we wish to assess, and $F$ will often be a nonparametric estimate $\hat{F}$ of the true sampling distribution $F_\tau$.

There are a number of important reasons why quadratic distances are central to the study of goodness of fit. These will be discussed in Section 2. A central goal here is to fill in some major gaps in the theory of quadratic distances. One major new result of this paper is the derivation of the limiting distribution for quadratic distances when used as goodness-of-fit tests in parametric models. These results depend on an appropriate spectral decomposition of the kernel $K$. We derive several new examples of such decompositions.

However, in many potential applications the numerical difficulty of determining the full spectral decomposition would make use of the limiting theory impractical. Our second set of major new results concerns the role of spectral degrees of freedom (sDOF), a concept introduced in this paper. We show that the limiting distributions involved are well approximated by chi-squared distributions when the degrees of freedom are large. Moreover, the sDOF are easily estimated empirically.

For kernel smoothing-based $L_2$ distances this is especially important because degrees of freedom are a more natural measure of the operating characteristics of the quadratic distance than are the bandwidth parameters. The literature on quadratic distances contains virtually no discussion of a concept we find critical. That is, in multivariate goodness of fit it is important to construct tuneable distances so that one can adjust the operating characteristics of the procedure to the dimension of the sample space and the sample size, much as one would do in a chi-squared analysis.

1.1. *The formal setup.* Let $\mathcal{S}$ be a sample space, with measurable sets $\mathcal{B}$, and let $du(s)$ be the canonical "uniform" measure on this space. The building block for our distance will be $K(s,t)$, a bounded, symmetric kernel function on $\mathcal{S} \times \mathcal{S}$. In analogy with matrix theory, a kernel is called *nonnegative definite* (NND), if the *quadratic form* $\iint K(s,t)\, d\sigma(s)\, d\sigma(t)$ is nonnegative for all bounded signed measures $\sigma$, and it is *conditionally* NND (i.e., CNND) if nonnegativity holds for all $\sigma$ satisfying the condition $\int d\sigma(s) = 0$.

Although our theoretical developments will be given for abstract spaces $\mathcal{S}$, it is important that for data calculations we will use discrete spaces. If $\sigma$ is finite discrete, with masses at $s_1, \ldots, s_m$, then the CNND requirements reduce to the conditional nonnegative definiteness of the matrix $\mathbb{K}$ having $i,j$ element $K(s_i, s_j)$.

DEFINITION 1. Given a CNND $K_G(s,t)$, possibly depending on $G$, the $K$-based quadratic distance between two probability measures $F$ and $G$ is

defined as

$$(1.1) \qquad d_K(F, G) = \iint K_G(s, t) \, d(F - G)(s) \, d(F - G)(t).$$

Note that the distance is well defined even when $F$ and $G$ do not have densities with respect to a common measure. The calculation for $d_K(F, G)$ can be written in the form

$$d_K(F, G) = K(F, F) - K(F, G) - K(G, F) + K(G, G),$$

where we have used the shortcut notation $K(A, B) = \iint K(s, t) \, dA(s) \, dB(t)$. We will call $d_K(\hat{F}, G)$ the *empirical distance* between the data and the $G$.

The discrete/matrix version of the problem will be of considerable statistical interest for its use in estimation. Let $F_\tau$ be the true sampling distribution and $\hat{F}$ the empirical distribution of a sample $x_1, \ldots, x_n$ from $F_\tau$. Let $\widetilde{\mathbb{K}_G}$ be the $n \times n$ *empirical representation* of the kernel $K_G$, having $ij$th element $K_G(x_i, x_j)$. In this case a quantity such as $\iint K_G(x, y) \, d\hat{F}(x) \, d\hat{F}(y) = \mathbf{1}^T \widetilde{\mathbb{K}_G} \mathbf{1} / n^2$ estimates $\iint K_G(x, y) \, dF_\tau(x) \, dF_\tau(y)$.

A possible practical limitation of quadratic distances is that numerical calculation of the distance requires twofold integration over the sample space. If the integrals are not explicit, one approach would be to perform Monte Carlo integration to calculate the distance, which in turn requires a simulation algorithm for the distributions involved. However, it is sometimes possible to choose a model-specific kernel that makes the distance calculation explicit and fast. This in turn enables one to construct test procedures that rely on other computationally intensive devices like bootstrapping.

**2. The central role of quadratic distance.** In this section we offer reasons that quadratic distance-based methods are central to goodness-of-fit inference.

2.1. *Important quadratic distances.* A number of classically important distances, such as Pearson's chi square or Cramér–von Mises, are quadratic distances. Other more recent examples can be found in Fan ([1997](#), [1998](#)), Fan, Zhang and Zhang ([2001](#)) and Zuo and He ([2006](#)).

*$L_2$ distances.* In $\mathcal{D} = \{1, 2, \ldots, N\}$ or $\mathcal{N} = \{0, 1, 2, \ldots\}$ one can use the "identity kernel" $K(s, t) = \mathbb{I}[s = t]$ and get the ordinary $L_2$ distance $\sum (f(i) - g(i))^2$. However, in $\mathcal{R}$ using the identity kernel gives the integral

$$\iint \mathbb{I}[x = y](f(x) - g(x))(f(y) - g(y)) \, dx \, dy,$$

which is identically zero. We will later show how to construct kernels that approximate the identity kernel, and hence the $L_2$ distance.

Similarly, the Pearson kernel

$$(2.1) \qquad K_G(s,t) = \frac{\mathbb{I}[s=t]}{\sqrt{g(s)g(t)}},$$

which nominally gives the Pearson distance $\int (f(s) - g(s))^2/g(s)\,du(s)$ between two densities, can be used in $\mathcal{D}$ or $\mathcal{N}$, but cannot be used in $\mathcal{R}$. This distance will be important to our story, as it is the quadratic distance approximant of Kullback–Leibler distance, as will be shown.

*An unconventional example.* A kernel that is quite popular as a smoothing kernel is the normal kernel with smoothing parameter $h$ [Silverman (1986)]. The special computational utility of this kernel derives from the convolution identity

$$(2.2) \qquad K_{h_1^2 + h_2^2}(x,y) = \int K_{h_1^2}(x,z) K_{h_2^2}(z,y)\,dz.$$

The identity implies if we use the normal kernel $K_{h^2}$ together with the normal model $G = K_{\sigma^2}(x,\mu)$ we obtain an explicit, no-integration-needed formula for the empirical distance given as $d_K(\hat{F}, G) = K_{h^2}(\hat{F}, \hat{F}) - 2n^{-1}\sum_i K_{h^2 + \sigma^2}(x_i, \mu) + K_{h^2 + 2\sigma^2}(\mu, \mu)$. This same computational facility carries over if $G$ is a finite mixture of normals.

2.2. *Relationship to $L_2$.* For a given symmetric kernel $K(s,t)$, there exists a symmetric kernel $K^{1/2}$ satisfying the relationship

$$\int K^{1/2}(s,r) K^{1/2}(r,t)\,du(r) = K(s,t).$$

(Existence will follow from the spectral decomposition that follows later.) For the normal kernel, (2.2) shows that $K_{h^2/2}$ is the square root kernel of the normal kernel $K_{h^2}$.

The square root operation leads us to a natural interpretation of the quadratic distance as an $L_2$ distance between smoothed densities.

PROPOSITION 1. *Let $K(s,t)$ be a symmetric, nonnegative definite kernel. Then*

$$d_K(F,G) = \int (f^*(z) - g^*(z))^2\,dz,$$

*where $f^*(z) = \int K^{1/2}(z,r)\,dF(r)$ and $g^*(z) = \int K^{1/2}(z,r)\,dG(r)$.*

PROOF. By reversal of order of integration. □

From the above proposition we see that if we use $\hat{F}$ instead of $F$, the empirical distance $d_K(\hat{F}, G)$ represents the $L_2(dz)$ distance between the kernel density estimator $f^*(z)$ and the smoothed $G$ distribution. Moreover, note that for the normal kernel the above relationship implies that the kernel is positive definite.

Conversely, any kernel smoothing problem can be turned into a quadratic distance problem. If, for example, $k_h(x - y)$ is a smoothing kernel on $\mathcal{R}$ that is used to construct the density estimator, the corresponding smoothed $L_2$ distance arises from the nonnegative definite kernel:

$$K(x, y) = \int k_h(x - z) k_h(z - y) \, dz.$$

This formula provides a simple way to generate CNND kernels from other kernels $k$.

2.3. *von Mises expansions.* We illustrate now that every smooth distance measure can be approximated, in a local sense, by a quadratic distance. To do this, we use the idea of von Mises expansion.

Consider the Kullback–Leibler distance $d(F, G) = \sum_{i=0}^{\infty} f(i) \ln[f(i)/g(i)]$ defined on $\mathcal{N}$. The influence function is

$$T'_{F^\circ}(s) = \ln f^\circ(s)/g(s) - \sum_{i=0}^{\infty} f^\circ(i) \ln[f^\circ(i)/g(i)].$$

Notice that the influence function is identically zero if "the null is true." Moreover, for the Kullback–Leibler distance, the Hessian is

$$T''_{F^\circ}(i, j) = \frac{\mathbb{I}[i = j]}{\sqrt{f^\circ(i) f^\circ(j)}}.$$

Thus, when the expansion point is $f^\circ = g$, the quadratic approximation to Kullback–Leibler is the Pearson chi-squared distance:

$$\sum_{i=0}^{\infty} f(i) \ln[f(i)/g(i)] \approx \sum_{i=0}^{\infty} \frac{[f(i) - g(i)]^2}{g(i)}.$$

**3. The decomposition theorem.** We now turn to discuss briefly the important role of spectral theory in determining the limiting distribution of the empirical quadratic distance $d_K(\hat{F}, G)$ between the data-based empirical distribution $\hat{F}$ and a hypothetical model $G$.

3.1. *Functional spectral decomposition.* Let $K(x, y)$ be a real-valued $\mathfrak{B}$-measurable positive definite kernel function on a measure space $(\mathcal{S}, \mathfrak{B}, M)$. The functional spectral decomposition of a kernel is similar to the spectral decomposition of a matrix with one very important exception: the functional spectral decomposition depends on the underlying measure $M$. In our usage, the distribution $M$ will usually be $F_\tau$, the true distribution of the data. If we are calculating the decomposition under the null hypothesis, $M$ will be $G$ for the simple hypothesis $H_0 : F_\tau = G$, and for composite null hypotheses $H_0 : F_\tau \in \{G_\theta\}$, $M$ will be one element of the parametric family of distributions. We will call $M$ the *baseline measure* and require that the kernel satisfies

$$(3.1) \qquad \int_{\mathcal{S}} \int_{\mathcal{S}} K(x,y)^2 \, dM(x) \, dM(y) < \infty.$$

This will hold for many typical examples because $M$ is a probability measure and $K$ is bounded. Such a kernel $K(x,y)$ generates a *Hilbert–Schmidt operator* on $L_2(M)$ through the *operation* $(Kg)(x) = \int K(x,y)g(y) \, dM(y)$. Our treatment here largely follows Yosida (1980).

THEOREM 3.1. *A nonnegative definite kernel $K$ satisfying* (3.1) *can be written as*

$$(3.2) \qquad K(x,y) = \sum_{j=1}^{\infty} \lambda_j \phi_j(x) \phi_j(y),$$

*where $\lambda_j$'s and $\phi_j(x)$'s are eigenvalues and corresponding normalized eigenvectors of $K$ under baseline measure $M$. The series in* (3.2) *converges strongly to $K$; that is, for every $g$ in $L_2$,*

$$\lim_{n \to \infty} \int_S \left( \int_S K(x,y)g(y) \, dM(y) \right.$$
$$\left. - \sum_{j=1}^{n} \int_S \lambda_j \phi_j(x) \phi_j(y) g(y) \, dM(y) \right)^2 dM(x) = 0.$$

*Moreover, $\lambda_j \geq 0$ since $K$ is NND.*

If $K(x,y)$ is real-valued and symmetric, then $K$ is a self-adjoint operator. The decomposition of $K(x,y)$ given in (3.2) corresponds to the spectral decomposition for a compact, self-adjoint operator, and will be called the $(K, M)$ spectral decomposition.

If $M$ equals the empirical measure $\hat{F}$, then $m(x_i) = 1/n$, the uniform density. Let $\widetilde{\mathbb{K}}$ be the $n \times n$ empirical matrix with $ij$th element $K(X_i, X_j)$. It is then clear that the $(K, \hat{F})$ eigendecomposition is just the same as the matrix eigendecomposition of the empirical kernel $\widetilde{\mathbb{K}}$ except that eigenvector normalization is changed from $\|\phi\|^2 = 1$ to

$$\int \phi^2(x) \, d\hat{F}(x) = n^{-1} \|\phi\|^2 = 1.$$

3.2. *Spectral trace of a kernel.* Fortunately, the most important attributes of the spectral decomposition can be calculated (or estimated) without obtaining the full decomposition. First, a consequence of the spectral decomposition theorem is that we can calculate the sum of squared eigenvalues by integration:

$$\sum_{j=1}^{\infty} \lambda_j^2 = \iint K(x,y)^2 \, dM(x) \, dM(y) < \infty.$$

We will denote the above quantity by $\text{trace}_M(K^2)$ due to its relationship to the matrix trace calculation.

The quantity $\text{trace}_M(K^2)$ is easily estimated from data. Suppose for measure $M$ we use the true distribution $F_\tau$, and that $X_1, \ldots, X_n$ is a sample from $F_\tau$. Then $\text{trace}_{F_\tau}(K^2)$ is estimated consistently by $\text{trace}_{\hat{F}}(K^2)$, which equals $\widetilde{\text{tr}(\mathbb{K}^2)}/n^2$, where in the last expression we have used tr to denote the standard matrix trace operation.

Many kernel functions satisfy even a stronger condition in that $\sum_{j=1}^{\infty} \lambda_j$, which we will write as $\text{trace}_M(K)$, is finite. The operators defined by those kernels are called nuclear. Under mild continuity assumptions one can also calculate $\text{trace}_M(K)$ without decomposition.

LEMMA 1. *Let $K(x,y)$ be a NND Hilbert–Schmidt kernel and let $\lambda_i$, $i = 1, 2, 3, \ldots$, denote the eigenvalues of the corresponding operator. Suppose that $K(x,y)$ is continuous at $(x,x)$ for almost all $x$ with respect to the measure $M$. Then, a necessary and sufficient condition for $\sum_{j=1}^{\infty} \lambda_j < \infty$ is that $\int K(x,x) \, dM(x)$ converges. Moreover, if $\sum_{j=1}^{\infty} \lambda_j < \infty$, then*

$$\sum_{j=1}^{\infty} \lambda_j = \int K(x,x) \, dM(x) = \text{trace}_M(K).$$

For the proof of this lemma see Yang ([2004](#)).

Once again, $\text{trace}_{F_\tau}(K)$ admits a simple consistent estimator, namely $\widetilde{\text{tr}(\mathbb{K})}/n$. These empirical estimators of trace quantities will be important later, as they enable one to approximate the limiting distributions of the test statistics through degrees-of-freedom calculations.

3.3. *An interpretation; plus centering.* Kernels and their representations are heavily used in support vector machines, where the eigenfunctions represent the "features" of importance in the problem, and the eigenvalues represent the weight attached to those features [Hastie, Tibshirani and Friedman ([2001](#))].

Similarly, in a statistical distance, we can write

$$d_K(F, G) = \iint K(s, t) \, d(F - G)(s) \, d(F - G)(t)$$

$$= \sum \lambda_i \left( \int \phi_j(s) \, dF(s) - \int \phi_j(s) \, dG(s) \right)^2,$$

so that the eigenvalues indicate the weight (importance) given the squared deviations in the *features*, which for us are the difference in expected values of the eigenfunctions under the two distributions.

However, there is an important detail missing. Because both $F$ and $G$ are probability measures, it is an easy exercise to show that $K(x, y)$ and $K^*(x, y) = K(x, y) + a(x) + a(y) + b$ both generate exactly the same quadratic distance $d_K(F, G)$, for any functions $a(x)$ and scalar $b$. However, $K$ and $K^*$ need not give the same spectral decomposition. Fortunately, statistical considerations point to a particularly natural choice for $a(x)$ and $b$ to use in the spectral decomposition. If $G$ is a hypothetical true model, then we should use the spectral decomposition of the following modified kernel.

DEFINITION 2. The *G-centered* kernel $K$, denoted by $K_{\text{cen}(G)}$, is defined as $K_{\text{cen}(G)}(x, y) = K(x, y) - K(x, G) - K(G, y) + K(G, G)$. When the identity of $G$ is clear from context, we will use notation $K_{\text{cen}}$, $K(x, G) = \int K(x, y) \, dG(y)$, and the terms $K(G, y)$ and $K(G, G)$ are similarly defined.

Note, by easy calculation, that

$$K_{\text{cen}}(x, G) = \int K_{\text{cen}}(x, y) \times 1 \, dG(y) = 0.$$

That is, the centering of $K$ has forced the function $\phi_1(x) = 1$ to be an eigenfunction of $K_{\text{cen}}$, with eigenvalue 0. As a consequence, by orthogonality to $\phi_1$, all the nonzero eigenfunctions have mean zero under $G$: $\int \phi_k(x) \, dG(x) = 0$.

The centering of the kernel is similar to a two-sided projection operation. If we are in $\mathcal{D}$, the discrete case, if $\mathbf{g}$ is the uniform density $\mathbf{1}/N$, as in the case of $\hat{F}$, and $\mathbf{1}\mathbf{1}^T/N$ is just the projection matrix $\mathbb{P}_{\mathbf{1}}$ that projects onto the space of constant vectors, then

$$(3.3) \qquad\qquad \mathbb{K}_{\text{cen}} = (\mathbb{I} - \mathbb{P}_{\mathbf{1}})\mathbb{K}(\mathbb{I} - \mathbb{P}_{\mathbf{1}}).$$

This "bilateral projection" formulation will later motivate the centering technique used when $G$ depends on estimated parameters.

In addition, if we wish to estimate nonparametrically the kernel after it has been centered by the true distribution $F_\tau$, we can empirically center the empirical kernel matrix $\widetilde{\mathbb{K}}$, obtaining

$$(3.4) \qquad\qquad \widetilde{\mathbb{K}_{\text{cen}}} = (\mathbb{I} - \mathbb{P}_{\mathbf{1}})\hat{\mathbb{K}}(\mathbb{I} - \mathbb{P}_{\mathbf{1}}).$$

We will later use this formula to estimate the total degrees of freedom.

3.4. *Examples of spectral decompositions.* In this section we will give several exact spectral decompositions.

3.4.1. *Poisson kernel.* In this subsection we construct a kernel by specifying the eigenvalues and eigenfunctions directly. The sample space will be the interval $[0, 2\pi)$ and the baseline measure $dM(x)$ will be the uniform probability density on this interval [i.e., $(2\pi)^{-1} dx$]. The eigenvalues for the kernel will have a geometrically decaying nature, $(\lambda_1, \lambda_2, \lambda_3, \ldots) = (1, \rho, \rho, \rho^2, \rho^2, \rho^3, \rho^3, \ldots)$, where $0 < \rho < 1$, with corresponding eigenfunctions $(1, \sqrt{2}\cos(x), \sqrt{2}\sin(x), \sqrt{2}\cos(2x), \sqrt{2}\sin(2x), \sqrt{2}\cos(3x), \sqrt{2}\sin(3x), \ldots)$. Written in terms of its spectral expansion, this gives the kernel

$$(3.5) \qquad K_\rho(\theta, \phi) = 1 + \sum_{k=1}^{\infty} 2\rho^k [\cos(k\theta)\cos(k\phi) + \sin(k\theta)\sin(k\phi)].$$

If one rewrites the cosine and sine terms in terms of complex exponential terms, one can use the geometric series formula to arrive at an explicit representation.

LEMMA 2.

$$K_\rho(\theta, \phi) = \frac{1 - \rho^2}{1 - 2\rho\cos(\theta - \phi) + \rho^2}$$

(3.6)

$$\text{where } 0 < \rho < 1 \text{ and } 0 \leq \theta, \phi < 2\pi.$$

Although not well known in statistics, this is the univariate version of the famous Poisson kernel. If we fix $\rho$ and $\phi$, then it becomes a density function in the variable $\theta$, with the parameter $\phi$ as a location parameter and $\rho$ as a dispersion parameter. This density has been used in statistics as a distribution on the unit circle, where it is known as the wrapped Cauchy distribution, first studied by Lévy (1939) and Wintner (1947).

In physics, it is the operator that gives the solution to the physical problem known as the "Dirichlet problem with boundary data" [e.g., Bhatia (2003)]. Additionally, it is a central tool in harmonic function theory [e.g., Axler, Bourdon and Ramey (2001)].

In this paper we will focus on the univariate version (3.6). Of importance to us here is that this distribution has a parameter, here $\rho$, that can be used to tune the degrees of freedom of the distance. Clearly, one could apply this kernel to distributions on any finite interval $[a, b]$ by a suitable location and scale change in the variables.

It is clear, using the infinite expansion (3.5) to do calculations, that after centering by the uniform distribution $M$, the Poisson kernel has the

decomposition

$$
\begin{aligned}
K_{\mathrm{cen}}(\theta, \phi) &= K(\theta, \phi) - 1 \\
&= \sum_{k=1}^{\infty} 2\rho^k [\cos(k\theta)\cos(k\phi) + \sin(k\theta)\sin(k\phi)].
\end{aligned}
\tag{3.7}
$$

It can also be shown that the appropriate convolution of two Poisson kernels is still a Poisson kernel, so this kernel is in many ways the natural analogue of the Gaussian one when one is considering data restricted to an interval.

3.4.2. *Normal kernel.* Of central importance to statistics is the spectral decomposition of the univariate normal kernel $K_{h^2}(x, y)$ when the baseline measure is $N(0, \sigma^2)$. A natural starting point is the Hermite polynomials. See Thangavelu ([1993]) for the relationships used here.

DEFINITION 3. The Hermite polynomials $H_n(x)$ are defined by the relationship

$$
H_n(x) = (-1)^n e^{x^2} \frac{d^n}{dx^n} e^{-x^2}.
$$

As candidates for the eigenfunctions we create a family of scaled and damped Hermite polynomials

$$
H_n(x; a, b) = H_n(ax) e^{-b^2 x^2/2},
$$

for $a$ and $b$ positive. These are useful because, using a classical identity for Hermite polynomials called *Mehler's formula*, we can create "spectral-like" expansions of $K_{h^2}$ in which the scaled and damped Hermite polynomials play the role of eigenfunctions (see the [Appendix] for the definition of $w^*$ and $\gamma^*$ used in the following results).

Let $\gamma_n^*(x)$ be $\gamma_n(x; a(w^*(r)), b(w^*(r)))$, as defined in ([A.1]).

THEOREM 3. *Under baseline measure $N(0, \sigma^2)$ the kernel $K_{h^2}(x, y)$ has the spectral decomposition $\sum_{n=0}^{\infty} \alpha\beta^n \gamma_n^*(x)\gamma_n^*(y)$, where $\beta = w^*(r)$ and*

$$
\alpha = \frac{(1 - w^{*2})^{1/2}}{2\sqrt{\pi} a(w^*)\sigma}.
$$

This representation shares with the Poisson kernel geometrically declining eigenvalues. It also captures the damped polynomial characteristic of the features used in the distance.

**4. Using distances for model assessment.** In this section we give some of the necessary theory behind testing-based model assessment.

4.1. *Estimation of the distance.* The next result gives a key property of the $G$-centered kernel.

PROPOSITION 2. *Let $F, G$ be two arbitrary distributions. Then the quadratic distance between $F, G$ can be written as*

$$d_K(F, G) = \iint K_{\text{cen}}(x, y) \, dF(x) \, dF(y).$$

This proposition shows that, for a fixed model $G$, the empirical distance $d_K(\hat{F}, G) = K_{\text{cen}(G)}(\hat{F}, \hat{F}) := V_n$ is a $V$-statistic [Serfling (1980)]. It can be calculated in matrix form as $\mathbf{1}^T \mathbb{K}_{\text{cen}} \mathbf{1}/n^2$. One can also unbiasedly estimate $d_K(F_\tau, G)$, where $F_\tau$ is the true distribution, by using the corresponding $U$-statistic:

$$(4.1) \qquad U_n = \frac{1}{n(n-1)} \sum_i \sum_{j \neq i} K_{\text{cen}}(x_i, x_j).$$

The fundamental distinction between $U_n$ and $V_n$ is the inclusion of the diagonal terms $K_{\text{cen}}(x_i, x_i)$, which have the nonzero expectation $\text{trace}_G(K_{\text{cen}})$. Under the null hypothesis $F_\tau = G$, the true distance $d(F_\tau, G)$ is zero, and $E_G(U_n) = 0$ but $E_G(V_n) = E[K_{\text{cen}}(X, X)]/n$, so that $\text{trace}_G(K_{\text{cen}})$ represents the biasing term.

4.2. *Under the null.* We start with the case where we have a prespecified null model $G$ that we wish to test, using as test statistic $V_n = d_K(\hat{F}, G)$ or the unbiased distance estimator $U_n(G)$. Letting $F_\tau$ denote the true distribution, the null hypothesis is $H_0 : F_\tau = G$. Given a spectral decomposition of the centered kernel $K_{\text{cen}}$ under $G$, say $K_{\text{cen}}(x, y) = \sum \lambda_i \phi_i(x) \phi_i(y)$, a heuristic derivation of the limiting distribution of $d_K(\hat{F}, G)$ is quite easy. Write

$$d_K(\hat{F}, G) = \iint \sum \lambda_i \phi_i(x) \phi_i(y) \, d\hat{F}(x) \, d\hat{F}(y)$$

$$= \sum_{i=1}^{\infty} \lambda_i (\bar{\phi}_i)^2,$$

where the $\bar{\phi}_i$ are averages of mean-zero, variance-1 variables that are uncorrelated over $i$. (Recall that the mean-zero property requires the use of the centered kernel.) The obvious conclusion is that

$$n V_n \xrightarrow{\text{dist}} \chi^*(\lambda), \qquad \lambda = (\lambda_1, \lambda_2, \ldots),$$

where $\chi^*(\lambda) = \sum \lambda_i Z_i^2$ is an infinite weighted sum of independent chi-squared variables. This is proved in Yang (2004).

The corresponding distributional result for the unbiased $U_n$ is that

$$\sqrt{n(n-1)}U_n \xrightarrow{\text{dist}} \chi^*_{\text{cen}}(\lambda),$$

where $\chi^*_{\text{cen}}(\lambda) = \sum \lambda_i(Z_i^2 - 1)$. This result was given formally for $U_n$ in Liu and Rao (1995), and holds under the condition that $\sum \lambda_i^2 = \text{trace}_G(K^2) < \infty$, which is weaker than the condition $\sum_{i=1}^{\infty} \lambda_i < \infty$ needed for $V_n$. Note that the result for $V_n$ cannot be improved upon because the distribution $\chi^*(\lambda)$ does not exist if $\sum \lambda_i = \infty$.

4.3. *Under composite nulls.* Next, consider the case where we wish to evaluate a parametric model $\{G_\theta : \theta \in \Omega\}$. We will assume that the elements $G_\theta$ of this model all have densities $g_\theta(x)$ with respect to a common measure $d\mu$. A natural test statistic for the validity of this model is $nV_n = nd_K(\hat{F}, G_{\hat{\theta}})$ (or the corresponding debiased statistic $U_n$) where $\hat{\theta}$ is a consistent estimator of $\theta$ under the null hypothesis $H_0 : F_\tau \in \{G_\theta\}$. If this method were applied to Pearson's kernel, for example, one would end up with Pearson's chi-squared test statistic.

The presence of the estimated parameter in $V_n$ necessarily makes finding the null distribution for general kernels $K$ more difficult, but we will show here that one can turn this problem into an eigendecomposition problem by artful centering of the kernel. Results similar to those presented here were derived by Fan (1998) for the special case of the weighted quadratic characteristic function distance.

Suppose that $p$-dimensional $\theta$ is being estimated by the maximum likelihood estimator $\hat{\theta}$. We will assume that it can be expressed as a solution to the set of $p$ likelihood equations

$$\sum \mathbf{u}(x_i; \theta) = 0,$$

where the likelihood scores $\mathbf{u}$ satisfy $E_\theta[\mathbf{u}(X; \theta)] = 0$. Notice that we are here using the maximum likelihood estimator for the problem, not the minimum quadratic distance estimator. The reason is that we anticipate that one would most likely use the quadratic distance fit assessment in conjunction with a maximum likelihood estimation procedure.

To find the distribution theory for the likelihood-estimated distance, we build a score-centered kernel from $K$ as follows. First, we construct the extended score vector $\mathbf{u}^* = (1, \mathbf{u}^T)^T$. We then define the extended information matrix for a single observation to be

$$\mathbb{J}^*_\theta = E_\theta[\mathbf{u}^*_\theta \mathbf{u}^{*T}_\theta].$$

We will then let $P^*$ be the kernel operator defined by

$$(4.2) \qquad P^*_\theta(x, y) = \mathbf{u}^*_\theta(x)^T \cdot \mathbb{J}^{*-1}_\theta \cdot \mathbf{u}^*_\theta(y).$$

The following formula shows that $P^*$ can be interpreted as the projection operator onto the extended space of likelihood scores:

$$(4.3) \qquad \int P_\theta^*(x, y) \mathbf{u}_\theta^{*T}(y) \, dG_\theta(y) = \mathbf{u}_\theta^{*T}(x).$$

The *score-centered kernel* is defined to be

$$
\begin{aligned}
(4.4) \quad K_{\text{scen}}^\theta &= (I - P_\theta^*) K (I - P_\theta^*) \\
&= K(x, y) - \int P_\theta^*(x, z) K(z, y) \, dG_\theta(z) \\
&\quad - \int K(x, z) P_\theta^*(z, y) \, dG_\theta(z) \\
&\quad + \iint P_\theta^*(x, z) K(z, w) P_\theta^*(w, y) \, dG_\theta(z) \, dG_\theta(w).
\end{aligned}
$$

The key feature of the score-centered kernel $K_{\text{scen}}^\theta$ is that it is $G_\theta$-orthogonal to the scores and the constant 1, as indicated next.

PROPOSITION 3. *The score-centered kernel satisfies*

$$\int K_{\text{scen}}^\theta(x, y) \mathbf{u}^*(y) \, dG_\theta(y) = \mathbf{0}.$$

This is easily proved using the definition (4.4) and repeated use of the projection property (4.3).

Note that the scores $\mathbf{u}$ are themselves orthogonal to the constant, that is, $\int \mathbf{u}(x) \cdot 1 \, dG_\theta(x) = \mathbf{0}$; therefore we could also have constructed the score-centered kernel by replacing (4.4) with

$$(4.5) \qquad (I - P_\theta) \cdot K_{\text{cen}(G_\theta)} \cdot (I - P_\theta),$$

where $P_\theta$ represents the projection onto the scores $\mathbf{u}$ instead of the extended scores $\mathbf{u}^*$.

In the discrete case, we can represent $P_\theta(i, j)$ by matrix $\mathbb{P}_\theta = \mathbf{u}_\theta \mathbb{J}_\theta^{-1} \mathbf{u}_\theta^T$, where $\mathbf{u}_\theta$ is the $N \times p$ matrix with entries $\partial_{\theta_j}[\log g_\theta(i)]$. We then get the matrix formula

$$\mathbb{K}_{\text{scen}}^\theta = (\mathbb{I} - \mathbb{P}_\theta \mathbb{D}_\theta) \cdot \mathbb{K}_{\text{cen}}^\theta \cdot (\mathbb{I} - \mathbb{D}_\theta \mathbb{P}_\theta),$$

where $\mathbb{D}_\theta$ is diagonal with diagonal entries $g_\theta(i)$.

The empirical distance between the data and the estimated model is then

$$d_K(\hat{F}, G_{\hat{\theta}}) = \iint K_{\text{scen}}^{\hat{\theta}}(x, y) \, d\hat{F}(x) \, d\hat{F}(y).$$

This can be verified by using the fact $\int \mathbf{u}_{\hat{\theta}}(x) \, d\hat{F}(x) = \mathbf{0}$ for maximum likelihood estimators.

This then leads to our main result.

THEOREM 4. *Given the regularity assumptions itemized in the proofs, under $G_\theta$ we have*

$$n\left[d_K(\hat{F}, G_{\hat{\theta}}) - \iint K_{\text{scen}}^\theta(x, y)\, d\hat{F}(x)\, d\hat{F}(y)\right] \overset{\text{prob}}{\longrightarrow} 0.$$

*If $K_{\text{scen}}^\theta$ has a spectral decomposition $\sum_i \lambda_i \phi_i(x)\phi_i(y)$ with finite trace, it follows that*

$$n d_K(\hat{F}, G_{\hat{\theta}}) \overset{\text{dist}}{\longrightarrow} \chi^*(\lambda),$$

*where $\lambda_1, \lambda_2, \ldots$ are the eigenvalues of the spectral decomposition of $K_{\text{scen}}^\theta$.*

In the Appendix we outline the steps in the proof. Notice that while there exists a corresponding U-statistic estimator of the distance, in the composite null hypothesis case it is no longer an unbiased estimator. One might still expect it to have slightly better operating characteristics.

**5. Spectral degrees of freedom.** We have now presented an asymptotic theory for quadratic distance methods that looks complicated and difficult to use. Except for certain carefully designed kernels, the spectral decomposition will be dependent on the underlying true model. It is likely there is not an explicit solution to the eigenequations. Even if the decomposition is known, the limiting distribution itself will depend on infinitely many $\lambda$ parameters.

These difficulties are not as severe as they first appear, because the key features of the spectral decomposition can be well summarized by the values of two scalar parameters called the *Pearson scale factor* and the *spectral degrees of freedom*. These two parameters are sufficient, in an asymptotic sense, for the description of the limiting distribution of the distance. As an additional bonus, they can easily be calculated for a model or estimated from the data without any spectral decomposition whatsoever.

5.1. *Pearson scaling and DOF.* Quadratic distances have no inherent scale. That is, replacing the kernel $K$ with $K^* = \alpha \cdot K$, for an arbitrary constant $\alpha$, creates a new distance that is equivalent to $K$ for most mathematical and statistical purposes.

Given a null measure $G$, we propose to rescale kernels so that they are as similar as possible to some standard kernel. The most natural standard kernel is the Pearson kernel. We will show that if we replace $K$ with $\alpha K$, where the scale factor is

$$(5.1) \qquad \alpha = \alpha_G(K) = \frac{\text{trace}_G(K)}{\text{trace}_G(K^2)} = \frac{\sum \lambda_i}{\sum \lambda_i^2},$$

then the quadratic distance generated by $\alpha K$ is scaled to match the Pearson kernel.

Given there is a fixed measure of interest $G$, we define a distance between the two kernels $K_1$ and $K_2$ via

$$\mathrm{trace}_G(K_1 - K_2)^2 = \iint (K_1(x,y) - K_2(x,y))^2 \, dG(x) \, dG(y).$$

Define the scaling factor $\alpha$ so that $\alpha K$ is as similar as possible to the Pearson kernel $Q$ by minimizing the distance

(5.2)
$$\begin{aligned}
&\mathrm{trace}_G(Q - \alpha K)^2 \\
&= \mathrm{trace}_G(Q^2) - 2\alpha(\mathrm{trace}_G(QK)) + \alpha^2 \, \mathrm{trace}_G(K^2).
\end{aligned}$$

Suppose for a moment we are in the finite discrete case, so we can write the Pearson kernel as $Q(x,y) = \mathbb{I}[x = y]/\sqrt{g(x)g(y)}$. In this case we have

$$\begin{aligned}
\mathrm{trace}_G(QK) &= \iint Q(x,y)K(x,y) \, dG(x) \, dG(y) \\
&= \int \frac{K(x,x)}{g(x)} g(x)g(x) \, du(x) \\
&= \mathrm{trace}_G(K).
\end{aligned}$$

Putting this into the expansion (5.2), the minimizing $\alpha$ is (5.1).

In other cases, if one minimizes the modified objective $\mathrm{trace}_G(-2\alpha KQ + \alpha^2 K^2)$, one again ends up with scale factor $\alpha_G(K)$.

5.2. *Spectral degrees of freedom.* We next define the *spectral degrees of freedom* (under $G$) of $K$ to be

(5.3)
$$\mathrm{DOF}_G(K) = \frac{\mathrm{trace}_G(K)^2}{\mathrm{trace}_G(K^2)} = \frac{(\sum \lambda_i)^2}{\sum \lambda_i^2}.$$

Note that $\mathrm{DOF}_G(K)$ equals $\mathrm{trace}_G(\alpha \cdot K)$, where $\alpha$ is the scale factor $\alpha_G(K)$. That is, the spectral degrees of freedom is just the sum of the eigenvalues of the rescaled kernel.

Fan, Zhang and Zhang (2001) found that the limiting normal distribution of their goodness-of-fit statistics had the mean-variance relationship of a scaled chi-squared distribution, and they used this to define the degrees of freedom of these tests. This relationship will be discussed in Section 5.4.

We will use a Satterthwaite approximation [Satterthwaite (1946)] to the $\chi^*(\lambda)$ distribution to interpret DOF. Recall that under the null empirical distance converges asymptotically in distribution to a linear combination of independent chi-squared random variables. Suppose we find scale $a$ and degrees of freedom DOF so that

$$E(\alpha d_K(\hat{F}, G)) = E(\chi^2_{\mathrm{DOF}}),$$

$$\mathrm{Var}(\alpha d_K(\hat{F}, G)) = \mathrm{Var}(\chi^2_{\mathrm{DOF}}).$$

Solving these two equations with respect to $a$ and DOF, we obtain

$$\alpha = \frac{2E(d_K(\hat{F}, G))}{\operatorname{Var}(d_K(\hat{F}, G))}$$

and

(5.4) $$\mathrm{DOF} = \frac{2E^2(d_K(\hat{F}, G))}{\operatorname{Var}(d_K(\hat{F}, G))}.$$

Using $E(d_K(\hat{F}, G)) = \operatorname{trace}_G(K) = \sum \lambda_i$ and $\operatorname{Var}(d_K(\hat{F}, G)) = 2\operatorname{trace}_G(K^2) = 2\sum \lambda_i^2$, we obtain that $\alpha$ is the Pearson scale factor and DOF is the same as defined in (5.3).

5.3. *Two examples.* In this subsection we will use two examples to illustrate calculation of the degrees of freedom. For point of comparison, we start with a classical quadratic distance, the Cramér–von Mises, which has a surprisingly small degrees of freedom. We then turn to the Poisson kernel as an example of the class of tuneable diffusion kernels. We show the degrees of freedom can be tuned to any value from 2 to infinity.

5.3.1. *Cramér–von Mises kernel.* The Cramér–von Mises kernel is given as $K(u, v) = 1 - \max(u, v)$ [Lindsay and Markatou (2002)]; its centered version is given as

$$K_{\mathrm{cen}}(u, v) = 1 - \max(u, v) - ((1 - u^2)/2) - ((1 - v^2)/2) + (1/3).$$

Using $G$ the uniform measure on $(0, 1)$, we obtain

$$\operatorname{trace}(K_{\mathrm{cen}}) = \int K_{\mathrm{cen}}(u, u)\, du = \int_0^1 (\tfrac{1}{3} + u^2 - u)\, du = \tfrac{1}{6}$$

and

$$\operatorname{trace}(K_{\mathrm{cen}}^2) = \int_0^1 \int_0^1 K_{\mathrm{cen}}^2(u, v)\, du\, dv = \tfrac{1}{90}.$$

Thus, the degrees of freedom for the centered Cramér–von Mises kernel are

$$\mathrm{DOF} = \frac{(1/6)^2}{(1/90)} = 2.5.$$

5.3.2. *The Poisson kernel.* For the Poisson kernel (3.5), let the baseline measure be uniform on $[0, 2\pi)$. Centering the kernel gives us

$$K_{\mathrm{cen}}(\theta, \phi) = K_\rho(\theta, \phi) - 1.$$

The following proposition gives the degrees of freedom of the centered Poisson kernel.

PROPOSITION 4. *The degrees of freedom of the centered Poisson kernel with respect to the uniform measure are given by*

$$\mathrm{DOF} = \frac{2(1+\rho)}{1-\rho}.$$

PROOF. From (3.7) the eigenvalues of the centered Poisson kernel with respect to the uniform measure are given by the set of functions $\{\rho, \rho, \rho^2, \rho^2, \ldots\}$. Now

$$\sum_{j=1}^{\infty} \lambda_i = 2\sum \rho^j = \frac{2\rho}{1-\rho}$$

and

$$\sum \lambda_i^2 = 2\sum \rho^{2j} = \frac{2\rho^2}{1-\rho^2}.$$

Therefore, the degrees of freedom are as given above.  □

When $\rho \to 0$ the degrees of freedom converge to 2, corresponding to the test that depends only on the first two eigenfunctions, whereas when $\rho \to 1$ the degrees of freedom diverge to infinity.

5.4. *Satterthwaite limit theory.* We now explore the relationship between the $\chi^*(\lambda)$ distribution, its Satterthwaite $\chi_K^2$ approximation and its normal approximation. A central assumption of this analysis is that we are considering kernels (like the normal or Poisson) with a tuning parameter $\eta$ such that the degrees of freedom become infinite as $\eta \to 0$.

The construction of $\chi^*(\lambda)$ as a sum of independent random variables suggests that we might hope for a central limit approximation for this distribution under the condition that the degrees of freedom are sent to infinity. If so, normality would imply that just two parameters would be sufficient to describe the distribution. We here give a simple sufficient condition for this result, and then go further. We will show that under the same conditions the Satterthwaite $\chi_{\mathrm{DOF}}^2$ approximation provides a two-parameter approximation that is always superior to the normal approximation.

We start by standardizing to mean zero and variance 1:

$$\chi_{\mathrm{std}}^*(\lambda) = \frac{\chi^*(\lambda) - \sum \lambda_i}{\sqrt{2\sum \lambda_i^2}}.$$

Since $\chi^*(\lambda)$ is a sum of independent variables, it is natural to study the cumulants of this distribution. Note that the cumulants of the standard normal, other than $r = 2$, are given by $\kappa_r(Z) = 0$, whereas $\kappa_2(Z) = 1$, so these are the cumulants we might hope to find in the limit.

To study the cumulants we define

$$\gamma_i = \frac{\lambda_i}{\sqrt{\sum \lambda_j^2}}.$$

Then we have $\sum \gamma_i^2 = 1$ and $\sum \gamma_i = \sqrt{\mathrm{DOF}(\lambda)}$. If we are considering the important special case of the $\chi_R^2$ distribution, then $R$ of the $\lambda_i$ are 1 and the rest are zero. This gives $\gamma_i = 1/\sqrt{R}$ for $R$ values of $i$ and 0 else.

The following lemma gives the cumulants for the standardized chi-star distribution. It arises from a straightforward calculation using the properties of cumulant-generating functions. Note that the cumulants of $\chi_1^2$ are given by $\kappa_r(\chi_1^2) = 2^{r-1}(r-1)!$.

LEMMA 5. *For $r \geq 2$ the $r$th cumulant of standardized $\chi^*(\lambda)$ is*

$$\kappa_r(\chi^*_{\mathrm{std}}(\lambda)) = \kappa_r(\chi_1^2) \cdot 2^{-r/2} \cdot \sum \gamma_i^r.$$

*For the $\chi_R^2$ distribution, an important special case, this gives $\kappa_r(\chi_{R,\mathrm{std}}^2) = \kappa_r(\chi_1^2) \cdot 2^{-r/2} \cdot R^{1-r/2}$.*

The degree of normality of the chi-star distribution can be measured by the departure of its cumulants from the normal values. In this case, we can show that the third cumulant (skewness) is the key factor.

LEMMA 6. *The normed cumulants*

$$\frac{\kappa_r(\chi^*_{\mathrm{std}}(\lambda))}{\kappa_r(\chi_1^2)2^{-r/2}} = \sum \gamma_i^r$$

*are decreasing in $r$ for $r = 2, 3, 4, \ldots$.*

PROOF. Each $\gamma_i$ is bounded above by 1, so $\gamma_i^r \geq \gamma_i^{r+1}$.  □

We have the following consequence of the last lemma: if we use a tuneable kernel with eigenvalues $\lambda_\eta$ depending on tuning parameter $\eta$, then all the cumulants of 3 and greater order converge to 0 as $\eta \to 0$ *if and only if* the skewness cumulant $\kappa_3(\chi^*_{\mathrm{std}}(\lambda_\eta))$ does. Indeed, if the skewness goes to zero, one can use the standard Taylor expansion proof to verify that

$$\chi^*_{\mathrm{std}}(\lambda_\eta) \xrightarrow{\mathrm{dist}} N(0,1) \qquad \text{as } \eta \to 0.$$

For the $\chi_R^2$ distribution, the skewness cumulant is $\kappa_3(\chi_{R,\mathrm{std}}^2) = 2^{3/2} \cdot R^{-1/2}$ reflecting its known convergence to normality. Below we will show that the skewness for the Poisson kernel goes to zero at the same rate in $R$, where $R$ is its spectral degrees of freedom.

If the zero-limit skewness property holds, one might ask whether there would sometimes be a preference for using the normal approximation over the Satterthwaite approximation to $\chi^*_{\text{std}}(\lambda)$. The answer is *never*, because the following lemma indicates that every cumulant of $\chi^*_{\text{std}}(\lambda)$ is closer to the Satterthwaite chi-squared cumulant than it is to the normal.

LEMMA 7. *Let $R$ be a positive integer. For $r \geq 3$, and for any $\chi^*(\lambda)$ distribution with $\text{DOF}(\lambda) = R$,*

$$\frac{\kappa_r(\chi^*_{\text{std}}(\lambda))}{\kappa_r(\chi^2_{\text{std},R})} \geq 1.$$

PROOF. See Appendix. □

That is, the cumulants are always larger in magnitude than the chi-squared ones, and further from zero, the normal theory value. This result also suggests that the magnitude of the skewness ratio

$$\frac{\kappa_3(\chi^*_{\text{std}}(\lambda))}{\kappa_3(\chi^2_{\text{std},R})}$$

could serve as a reasonable index of the relative chi-squaredness of the chi-star distribution when the degrees of freedom are large. Additionally, the limit of this ratio as $R$ becomes infinite could serve as a single number summary.

For the Poisson kernel example, with $\rho = e^{-\eta}$ and using the uniform baseline density, we get

$$\frac{\kappa_3(\chi^*_{\text{std}}(\lambda))}{\kappa_3(\chi^2_{\text{std},R})} = \frac{(1 + e^{-\eta})^2}{(1 + e^{-\eta} + e^{-2\eta})},$$

a term which converges to $4/3$ as $\eta \to 0$. That is (and we found this surprising), with geometrically declining eigenvalues the skewness of standardized $\chi^*(\lambda)$ lies closer to the chi-square's skewness than the latter's does to the normal value of 0. In general, for the Poisson kernel the ratio of $r$th cumulants converges to $2^{r-1}/r$, showing that the $r$th cumulant is the same magnitude as the chi-square: $O(R^{1-r/2})$, where $R$ is the degrees of freedom.

6. **Final comments.** Quadratic distances with tuning parameters are in many ways like smooth chi-squared goodness-of-fit tests: the $L_2$ relationship suggests that, as an alternative to constructing a finite set of bins, we are creating an infinite number and averaging across their deviations. The spectral degrees of freedom concept is meant to be a tool to help statisticians exploit this analogy.

If we accept this analogy, then the choice of the degrees of freedom in a problem is like the choice of the number of cells in the chi-squared test: it is clearly extremely important in determining the power of the test against interesting alternatives, but it is also very hard to devise hard and fast rules about its choice. That is because the nature of the interesting alternatives may not be clear to the user. We might add, that provided one is using the test statistics as an exploratory tool, there is no reason one would not consider a range of interesting degrees of freedom as a means of exploring the possible deviations from the model at various scales of smoothing. (We think that informal/exploratory model confirmation is widely used, and this would simply be another instance.)

What then is an interesting range for degrees of freedom? At this time, we can only offer a heuristic analysis based on chi-squared tests. In a very general sense, increasing the number of cells in such a test, and therefore the degrees of freedom, will create a gain in sensitivity to deviations that are localized within a single small area (like a bump in the density), but also create increased variability of the test statistic that causes it to lose power against more global alternatives that create a small shift in probability in many cells.

In a chi-squared test one would want, even if searching for small local deviations, some minimum sample counts per cell in order to cut variability. If that minimum were 5, one would never have more than $n/5$ degrees of freedom. This is a number we have used as a rough upper bound when we investigated a problem.

On the other hand, just as a chi-squared test with two cells would be too coarse for most purposes, one should avoid having too small a degrees of freedom. In this regard, the dimension of the sample space is important. To illustrate, if one were to provide a one-degree-of-freedom test on each marginal distribution in a $D$-dimensional data set, then one has used $D$ degrees of freedom. To also test all the bivariate marginals would take an additional $\binom{D}{2}$ degrees of freedom. Based on this heuristic, we have used $\binom{D+1}{2}$ as a very rough lower bound when investigating multivariate data sets.

One important issue we have not touched upon in this paper is that of power. It is very difficult to draw broad conclusions about test procedures based on their power characteristics because the dimension of the alternative space is infinite, and it is inevitable that the identity of the best performing test will be highly dependent on the alternative that is chosen. Spitzner (2006) developed a detailed simulation study that compared a variety of testing strategies for combining quadratic tests into a single test statistic. Best power? The answer depended on the structure of the alternative; Fan's adaptive Neyman strategy [Fan (1996)] worked the best overall in Spitzner's particular simulation settings, but was not a universal winner.

## APPENDIX: PROOFS AND LIMITING DISTRIBUTIONS

### A.1. Proofs for Theorem 3.

LEMMA A.1. *For any $w \in (0,1)$, let $a = a(w) = \sqrt{\frac{(1-w)^2}{2h^2 w}}$ and $b = b(w) = \sqrt{\frac{1-w}{h^2}}$. Then*

$$K_{h^2}(x,y) = \frac{1}{\sqrt{2\pi} h} \sum_{i=0}^{\infty} \left[ \frac{w^n (1-w^2)^{1/2}}{2^n n!} \right] H_n(x;a,b) \cdot H_n(y;a,b).$$

PROOF. Mehler's formula states that for $w \in (0,1)$

$$\sum_{n=0}^{\infty} c_n(w) H_n(x) H_n(y) = \exp\left( \frac{2xyw - (x^2 + y^2)w^2}{1 - w^2} \right),$$

where

$$c_n(w) = \left[ \frac{w^n (1-w^2)^{1/2}}{2^n n!} \right].$$

A series of algebraic manipulations leads to the desired representation. □

The above formula is not a spectral representation unless we can choose $a$ and $b$ so that the $H_n(x;a,b)$ terms are orthogonal under the normal measure. The following gives us the necessary condition on $a$ and $b$.

LEMMA A.2. *Let $a$ and $b$ be two positive scalars satisfying $a^2 - b^2 = (2\sigma^2)^{-1}$. The functions $\gamma_0, \ldots, \gamma_n, \ldots$ defined by*

$$(A.1) \qquad \gamma_n(x;a,b) = H_n(x;a,b) \sqrt{\frac{a\sigma}{2^{n-1/2} n!}}$$

*are orthonormal under the measure $M = N(0, \sigma^2)$.*

PROOF. We start with the fundamental identity

$$\int_{\mathcal{R}} H_m(x) H_n(x) e^{-x^2} dx = \mathbb{I}[m = n] 2^n n! \sqrt{\pi}.$$

With a change of variables $x = ay$, the left-hand side becomes

$$LHS = a \int_{\mathcal{R}} H_m(ay) H_n(ay) e^{-a^2 y^2} dy$$

$$= a \int_{\mathcal{R}} H_m(ay) e^{-b^2 y^2/2} H_n(ay) e^{-b^2 y^2/2} e^{-(a^2 - b^2)y^2} dy$$

$$= a \int_{\mathcal{R}} H_m(y;a,b) H_n(y;a,b) e^{-y^2/2\sigma^2} dy.$$

We therefore have

$$\frac{1}{\sqrt{2\pi}\sigma} \int_{\mathcal{R}} H_m(y;a,b) H_n(y;a,b) e^{-y^2/2\sigma^2} \, dy = \frac{\mathbb{I}[m=n] 2^n n! \sqrt{\pi}}{\sqrt{2\pi}\sigma a},$$

as needed. $\square$

The final trick is to try to select the scalar $w^*$ in the first lemma such that the functions $a(w^*)$ and $b(w^*)$ that are defined there satisfy the condition $a^2 - b^2 = (2\sigma^2)^{-1}$ of the second lemma. Let $r = h^2/\sigma^2$, the ratio of the kernel and baseline variances.

LEMMA A.3. *Set $w^*(r) = 1 - \frac{1}{2}[\sqrt{4r + r^2} - r]$. Then $w^*(r)$ decreases monotonely from 1 to 0 as a function of $r$, for $r \in (0, \infty)$. For any $h^2$ and $\sigma^2$, we have $a(w^*(r))^2 - b(w^*(r))^2 = (2\sigma^2)^{-1}$.*

PROOF. The function $w^*(r)$ is the left-hand root of the quadratic equation $rw = (1 - w)^2$. Inspecting the plot of the two sides of this quadratic equation verifies the listed functional properties. The last equality is easy algebra. $\square$

**A.2. Proofs for Theorem 4.** First, we show that score-centering implies mean-centering of the derivatives of the kernel.

PROPOSITION A.1. *If the kernel is score-centered under $F_\tau = G_\theta$, then under regularity conditions*

$$(A.2) \qquad \int (\nabla_\theta K_{sc}^\theta(x,y)) \, dG_\theta(y) = \mathbf{0}.$$

*In addition,*

$$(A.3) \qquad \int\int \nabla^2 K_{sc}^\theta(x,y) \, dG_\theta(x) \, dG_\theta(y) = \mathbf{0}.$$

The proof can be easily obtained by differentiation under the integral sign. These mean-zero properties can then be used to show the following:

PROPOSITION A.2. *Under regularity conditions found in the proof,*

$$n \int\int [K_{sc}^{\hat{\theta}}(x,y) - K_{sc}^\theta(x,y)] \, d\hat{F}(x) \, d\hat{F}(y) \xrightarrow{\text{prob}} 0.$$

PROOF. We plug the following Taylor expansion:

$$\begin{aligned}
K_{sc}^{\hat{\theta}}&(x,y) - K_{sc}^\theta(x,y) \\
&= (\hat{\theta} - \theta)^T [\nabla K_{sc}(x,y)] \\
&\qquad + \tfrac{1}{2}(\hat{\theta} - \theta)^T [\nabla^2 K_{sc}^\theta(x,y)](\hat{\theta} - \theta) + (rem)
\end{aligned}$$

into the above expression. We then assume root-$n$ consistency of the MLE, so that $\sqrt{n}(\hat{\theta} - \theta)$ converges in distribution. We assume that if $X$ and $Y$ are independent from $G_\theta$, the kernels $\nabla K_{\mathrm{sc}}(X, Y)$ and $\nabla^2 K_{\mathrm{sc}}(X, Y)$, which are mean zero from Proposition A.1, have finite variance, as do $\nabla K_{\mathrm{sc}}(X, X)$ and $\nabla^2 K_{\mathrm{sc}}(X, X)$. If so, then $n \iint \nabla K_{\mathrm{sc}}(x, y) \, d\hat{F}(x) \, d\hat{F}(y)$ and $n \iint \nabla^2 K_{\mathrm{sc}}(x, y) \, d\hat{F}(x) \, d\hat{F}(y)$ converge in distribution. This then assures that the first two terms are of stochastic order $O_p(n^{-1/2})$ and $O_p(n^{-1})$, respectively. The remainder term is then no larger than $O_p(n^{-1/2})$ under the assumption that the elements of the arrays $\nabla^3 K_{\mathrm{sc}}^{\theta^*}(X, Y)$ and $\nabla^3 K_{\mathrm{sc}}^{\theta^*}(X, X)$ are bounded by finitely integrable functions for $\theta^*$ in a neighborhood of $\theta$. $\square$

### A.3. Proofs for Lemma 7.

We start by proving the result when the eigenvalue sequence is finite in length, say $\gamma_1, \ldots, \gamma_N$. We can then write $\sum_{i=1}^N \gamma_i^r / N = \sum_{m=1}^M \pi_m \cdot a_m^r$, where $a_1, \ldots, a_M$ represent the $M$ distinct values possible among the $\gamma_i$ and the $\pi_m$ represent the counts for each $a_m$, divided by $N$. The expression $\sum \pi_m a_m^r$ is therefore the $r$th moment of the distribution that puts mass $\pi_m$ at support point $a_m$. For this distribution we know the first two moments: $\sum \pi_m a_m^1 = \sqrt{R}/N$ and $\sum \pi_m a_m^2 = 1/N$. We wish to know the minimum possible value of the $r$th moment over the possible distributions represented by $\pi_m$ and $a_m$.

We enlarge the class of allowable distributions to include every distribution $P$ with its support in $[0, 1]$. Under the theory of moments, the solution to this optimization problem is an extremal distribution having index $3/2$. That is, the optimal $P$ has two support points, one of which is 0 or 1. We can exclude 1 because this would maximize the $r$th moment, leaving us with one support point of 0. However, the $\chi^2_{\mathrm{std}, R}$ distribution has the eigenvalue distribution $P$ of index $3/2$, putting all its probability on the two support points 0 and $1/\sqrt{R}$ with masses $(N - R)/N$ and $R/N$, respectively, and so it has the extremal $r$th moment. The theorem is concluded by taking limits as $N \to \infty$.

B. G. LINDSAY
DEPARTMENT OF STATISTICS
PENNSYLVANIA STATE UNIVERSITY
UNIVERSITY PARK, PENNSYLVANIA 16802
USA
E-MAIL: bgl@psu.edu

S. RAY
DEPARTMENT OF MATHEMATICS
  AND STATISTICS
BOSTON UNIVERSITY
111 CUMMINGTON STREET
BOSTON, MASSACHUSETTS 02215
USA
E-MAIL: sray@math.bu.edu

M. MARKATOU
DEPARTMENT OF BIOSTATISTICS
722 WEST 168TH STREET
ROOM 632
COLUMBIA UNIVERSITY
NEW YORK, NEW YORK 10032
USA
E-MAIL: mm168@columbia.edu

K. YANG
CYTOKINETICS, INC.
280 EAST GRAND AVENUE
SOUTH SAN FRANCISCO, CALIFORNIA 94080
USA
E-MAIL: keyang@gmail.com

S.-C. Chen
439 PSA
Department of Mathematics
    and Statistics
Arizona State University
Tempe, Arizona 85287
USA
and
Department of Statistics
National Cheng Kung University
Taiwan
ROC
E-mail: scchen@math.asu.edu